\documentclass[12pt]{article}
\usepackage{amsmath,amssymb,amsfonts,euscript}

\newtheorem{Theorem}{Theorem}[section]
\newtheorem{Lemma}[Theorem]{Lemma}
\newtheorem{Corollary}[Theorem]{Corollary}
\newtheorem{Proposition}[Theorem]{Proposition}
\newtheorem{Remark}[Theorem]{Remark}
\newtheorem{Example}[Theorem]{Example}
\newtheorem{Examples}[Theorem]{Examples}
\newtheorem{Conjecture}[Theorem]{Conjecture}
\newtheorem{Definition}[Theorem]{Definition}
\newtheorem{Question}[Theorem]{Question}
\newcommand{\rar}{\rightarrow}

\newcommand{\surjects}{\twoheadrightarrow}

\def\sqr#1#2{{\vcenter{\hrule height.#2pt
        \hbox{\vrule width.#2pt height#1pt \kern#1pt
            \vrule width.#2pt}
        \hrule height.#2pt}}}
\def\demo{\noindent{\bf Proof. }}
\def\square{\mathchoice\sqr64\sqr64\sqr{4}3\sqr{3}3}
\def\qed{\hspace*{\fill} $\square$}

\def\Rees{{\cal R}}
\def\hht{{\rm ht}\,}
\def\spec{{\rm spec}\,}
\def\XX{{{\mathbb X}}}

\def\fm{{\mathfrak m}}
\def\NN{\mathbb N}
\def\XX{{{\bf X}}}
\def\cl#1{{\cal #1}}

\title{Aluffi Torsion-Free Ideals}
\author{Abbas  Nasrollah Nejad \\   Rashid Zaare-Nahandi \\
Institute for Advanced Studies in Basic Sciences (IASBS)\\
P.O.Box 45195-1159,
Zanjan 45137-66731, Iran\\
abbasnn@gmail.com, rashidzn@iasbs.ac.ir}

\date{}

\begin{document}

\maketitle


\begin{abstract}
A special class of algebras which are intermediate between the symmetric and the Rees algebras of an ideal was introduced by P. Aluffi in 2004 to define characteristic cycle of a hypersurface parallel to conormal cycle in intersection theory. These algebras are recently investigated by A. Nasrollah Nejad and A. Simis who named them Aluffi algebras. For a pair of ideals $J\subseteq I$ of a commutative ring $R$, the Aluffi algebra of $I/J$ is called Aluffi torsion-free if it is isomorphic to the Rees algebra of $I/J$.
In this paper, ideals generated by 2-minors of a $2\times n$ matrix of linear forms and also edge ideals of graphs are considered and some conditions are presented which are equivalent to Aluffi torsion-free property of them. Also many other examples and further questions are presented.
\end{abstract}

\section*{Introduction}
In the remarkable paper~\cite{aluffi},  Paolo Aluffi introduced  an intermediate graded algebra
between a symmetric algebra and the  Rees algebra which he called quasi symmetric algebra. His purpose was to describe the characteristic cycle of a hypersurface, parallel to well known conormal cycle in intersection theory.
A. Nasrollah Nejad and A. Simis in~\cite{AA}
called such an algebra the \textit{Aluffi algebra}. Given a commutative ring $R$ and
ideals $J\subset I\subset R$, the Aluffi algebra of $I/J$ is defined by
$${\cal A}_{R/J}(I/J):={\cal S}_{R/J}(I/J)
\otimes_{{\cal S}_R(I)}\Rees_R(I).$$

The Aluffi algebra is squeezed as ${\cal S}_{R/J}(I/J)\surjects {\cal A}_{{R/J}}(I/J)\surjects\Rees_{R/J}(I/J)$ and is
moreover a residue ring of the ambient Rees algebra $\Rees_R(I)$. The kernel of the right hand surjection is called the module of Valabrega-Valla as defined in \cite{VaVa} which is the torsion of the Aluffi algebra \cite{AA}. Thus, provided that $I$ has a regular element modulo $J$, the Rees algebra of $I/J$ is the Aluffi algebra modulo its torsion. The question which motivated this paper is: when is the surjection ${\cal A}_{{R/J}}(I/J)\surjects\Rees_{R/J}(I/J)$ an isomorphism. For importance of this question in commutative algebra and intersection theory, we call a pair of ideals $J\subset I$, Aluffi torsion-free
 if the surjection ${\cal A}_{{R/J}}(I/J)\surjects\Rees_{R/J}(I/J)$ is injective.

Some important examples of Aluffi torsion-free pairs have been appeared explicitly in the following
two results. The first one is due to Huneke \cite{Huneke} who states an ideal $I$ which its extension $(I+J)/J$ on the quotient ring $R/J$ is generated by a $d$-sequence.
The second one is due to Herzog, Simis and Vasconcelos and is what they called "Artin-Rees lemma on the nose" \cite{Trento}. They have considered that, both ideals $I$ and $I/J$ are of linear type over $R$ and $R/J$, respectively. By the structure of the Aluffi algebra, it is shown in \cite{AA} that the assumption in the second result to the effect that $I$ be of linear type over $R$ does not intervene the result. Nasrollah Nejad and Simis in \cite{AA} give necessary and sufficient conditions for these algebras to be isomorphic in terms of
$I$-standard basis of $J$ and also relates this isomorphism with the relation type number
of $I/J$ over $R/J$ and the Artin-Rees number of $J$ relative to $I$.

In geometric settings, let $ X\stackrel{i}\hookrightarrow Y \stackrel{j}\hookrightarrow Z$  be closed embeddings of schemes with $J\subset I\subset R$ the ideal sheaves of $Y$ and $X$ in $Z$, respectively. Let $\widetilde{Z}={\rm Proj} (\Rees_R(I))\stackrel{\pi}\rar Z$ be the blowup of $Z$ along $X$ and $\widetilde{Y}={\rm Proj} (\Rees_{R/J}(I/J))$ be the blowup of $Y$ along $X$. Note that $\widetilde{Y}$ embeds in $\widetilde{Z}$ as the strict transform of $Y$ under $\widetilde{Z}\stackrel{\pi}\rar Z$.
Let $E=\pi^{-1}(X)$ be the exceptional divisor of the blowup. Then, $E$ is a subscheme of $\pi^{-1}(Y)$. Let $\mathfrak{R} = \mathfrak{R}(E, \pi^{-1}(Y))$ be the residual scheme of $E$ in $\pi^{-1}(Y)$. Here "residual" is taken in the sense of \cite[Definition 9.2.1]{fulton}. In terms of the ideal sheaves, $\mathfrak{R}$ is characterized by the equation $\mathcal{I}_{\mathfrak{R}}.\mathcal{I}_{E}=\mathcal{I}_{\pi^{-1}(Y)}$, where ${\mathcal I}_E, {\mathcal I}_{\pi^{-1}(Y)}$ are respectively the ideals of $E$ and $\pi^{-1}(Y)$ in $\tilde{Z}$. Aluffi in  \cite[Throrem 2.12]{aluffi} proved that ${\rm Proj}({\cal A}_{{R/J}}(I/J))=\mathfrak{R}(E, \pi^{-1}(Y))$. Fulton in \cite[B. 6.10]{fulton} shows that if $i$ and $j$ are regular embeddings, then $\mathfrak{R}= \widetilde{Y}$ which is equivalent to say that $J\cap I^{n}=JI^{n-1}$ for all sufficiently large $n$. S. Keel in \cite[Theorem 1]{keel} shows that this result holds as long as $X\hookrightarrow Y$ is a linear embedding and $Y\hookrightarrow Z$ is a regular embedding.
The goal of the present work is to find some examples of Aluffi torsion-free pairs which are in the main streams of researches in Commutative Algebra and Algebraic Geometry. We classify completely the ideals generated by 2-minors of a $2\times n$ matrix of linear forms and edge ideals of graphs in terms of the Aluffi torsion-free property.

In Section 2, we consider $J$ as an ideal generated by 2-minors of a $2\times n$ matrix of linear forms and $I$ stands for the Jacobian ideal of $J$. We prove that the pair $J\subseteq I$ is Aluffi torsion-free if and only if in the Kronecker-Weierstrass normal form of the matrix, there is no any Jordan block. More precisely, Theorem~\ref{matrix} asserts that these conditions are equivalent to say that $I_r(\Theta) = \fm^r$, where $r$ is codimension of $J$, $\Theta$ stands for the Jacobian matrix of $J$ and $\fm$ is the homogeneous maximal ideal of $k[\XX]=k[x_1,\ldots,x_n]$. This motivates us to conjecture that, if $J\subset k[\XX]$ is an ideal of codimension $r\geq 2$, generated by 2-forms, and if $I$ denotes the ideal generated by $r$-minors of the Jacobian matrix $\Theta$ of $J$, then $I$ is $\fm$-primary if and only if $I=\fm^r$ (Conjecture~\ref{conj1}).

Section 3 is devoted to find conditions for edge ideal of a graph  and its Jacobian ideal to be Aluffi torsion-free pair. In this regard, we give some necessary and sufficient conditions for graphs equivalent to the Aluffi torsion-free property. Finally, we present several examples of graphs
which are Aluffi torsion-free or not.

In the last section, some more examples of important ideals are considered and some questions for further steps are posed.

Some of the results of this paper have been conjectured after explicit computations performed by the computer algebra systems Singular~\cite{Sing} and CoCoA~\cite{coco}.

\section{Torsion-Free Aluffi algebras}
 Let $R$ be a commutative ring and $I$ an ideal of $R$.
The two most common and important commutative algebras related to the ideal $I$ are the
symmetric algebra ${\cl S}_R(I)$ and the Rees algebra ${\cl R}_R(I)$.
Recall that these algebras are defined as
$${\cl R}_R(I):=\bigoplus_{t\geq 0} I^tu^t\simeq R[Iu]\subset R[u],\ \ \ {\cl S}_R(I) := \bigoplus_{t\geq 0}{\cl S}_R^t(I),$$
where, ${\cl S}_R^t(I)={T^t_R(I)}/{(( x\otimes y - y\otimes x) \cap T^t_R(I))}$ and $T^t_R(I)$ is the tensor algebra of order $t$.
The definition of ${\cl R}_R(I)$ immediately implies that, it is torsion-free over the base ring $R$.
A natural surjection of standard $R$-graded algebras arises from the definition:
\begin{equation}\label{Sym_to_Rees}
{\cl S}_R(I)\surjects {\cl R}_R(I).
\end{equation}
This map is injective locally on the primes $p\in \spec(R)$ such that $I\not\subseteq p$.
It follows from the general arguments that, provided that $I$ has some regular elements, the kernel
is the $R$-torsion submodule (ideal) of the symmetric algebra.
If the map in (\ref{Sym_to_Rees}) is injective, one says that the ideal $I$ is of linear type,
a rather non-negligible notion in parts of  syzygy theory of ideals.

\begin{Definition}\rm (\cite{AA})
Let $R$ be Notherian  and $J\subset I$  be ideals of $R$. The Aluffi algebra of $I/J$ is
$${\cl A}_{{R/J}}(I/J) := {\cl S}_{R/J}(I/J)
\otimes_{{\cl S}_R(I)}
\Rees_R(I).$$
 \end{Definition}
We have the following  surjections:
$${\cl S}_{R/J}(I/J)\surjects {\cal A}_{{R/J}}(I/J)\surjects\Rees_{R/J}(I/J). $$
The kernel of the second
surjection is the so-called module of Valabrega-Valla (see
\cite{VaVa}, also \cite[5.1]{Wolmbook1}) which is:
\begin{equation}\label{vava}
 {\cl V}\kern-5pt {\cl V}_{J\subset I}=\bigoplus_{t\geq 2} \frac{J\cap
I^t}{JI^{t-1}}.
\end{equation}
Of course, as an ideal, this kernel is generated by finitely many homogeneous
elements, but as a graded $R/J$-module, it is conceivable that it may fail this property. By \cite[Proposition 2.5]{AA} the Valabrega-Valla module gives the torsion of the Aluffi algebra.
\begin{Definition}\label{torsionfree}\rm
A pair of ideals $J\subset I$ of a ring $R$ is said to be Aluffi torsion-free
if the map ${\cal A}_{R/J}(I/J)\surjects\Rees_{R/J}(I/J)$ is injective.
\end{Definition}

Note that by \cite[B. 6.10]{fulton} and (\ref{vava}), a pair of ideals $J\subset I$ is Aluffi torsion-free if and only if $J\cap I^{n}=JI^{n-1}$ for all positive integers $n$.
\begin{Example}\rm
Let $J\subset I$ be proper ideals of a ring $R$.
If $I/J$ is of linear type over $R/J$ (e.g. if $I$ is generated by
a regular sequence modulo $J$), then $J\subset I$ is Aluffi torsion-free.
\end{Example}

\begin{Example}\rm
Let $a_1,\ldots,a_r$ be a regular sequence in a Noetherian ring $R$ and let $I=\langle a_1,\ldots,a_r\rangle$.
Then, for each $i=1,\ldots,r$, the pair $J=(a_1^n,\ldots,a_i^n)\subset I^n$  is  Aluffi torsion-free.

To see this, by induction on $t$ we show that $J\cap I^{nt}=JI^{n(t-1)}$ for every $t\geq 1$.
For $t=1$ the conclusion is obvious.
Let $t>1$, then by the inductive assumption
$$J \cap I^{nt}=I^{nt}\cap J \cap I^{n(t-1)}=I^{nt}\cap JI^{n(t-2)}=I^{nt}\cap \Delta,$$
where $\Delta$ is the ideal generated by the elements $a_1^{s_1}\ldots a_r^{s_r}$
such that $s_1+\ldots+s_r=n(t-1)$ and $s_j\geq n$ for $1\leq j\leq i$. Since $a_1,\ldots,a_r$
is a regular sequence, if $f(x_1,\ldots,x_i)$ is a homogenous polynomial of degree $n(t-1)$
over $R$ such that $f(a_1,\ldots,a_i)\in I^{nt}$, then, all coefficients of $f$ must be in $I$.
Therefore,
$$ J\cap (I^n)^t=J\cap I^{nt}=I^{nt}\cap \Delta=I\Delta =JI^{n(t-1)}.$$
\end{Example}

\begin{Lemma}\label{forms}
Let $R=k[\XX]$ and $J\subset R$  be an ideal  generated by forms of the same degree $d\geq 1$.
Then, $J\cap \fm^{rt}\subset
J\fm^{r(t-1)}$ for every $t\geq 0$ and  $r\geq d$.
\end{Lemma}
\demo
Let $f_1,\ldots,f_m$ be generators of $J$ and let $F$ be a form on $f_i$'s such that $F\in \fm^{rt}$.
Then $F=\sum _{i=1}^{m}g_if_i$, where $g_i=\sum{a_{\alpha}}\XX^{\alpha}\in R_{rt-d+\delta}$
for $\delta\geq 0$. Since $R_{rt-d+\delta}=R_{r-d+\delta}.R_{rt-r}$, we can rewrite $g_i$ as
$$g_i=\displaystyle\sum_{{|\alpha|=r-d+\delta}\atop {|{\beta}|=rt-r}} \
 a_{\alpha,\beta}\ \XX^{\alpha+\beta},\quad {\rm hence}\quad
F= \sum_{{|{\alpha}|=r-d+\delta}\atop\ }\ \XX^{{\alpha}}
\left(\sum_{i=1\atop |{\beta}|=rt-r}^{s}(\XX^{{\beta}})f_i\right).$$
Therefore, $F\in J\fm^{rt-r}$, as required.
\qed

\medskip

Let $R=k[\mathbf{X}]$ be the $\NN$-graded polynomial ring over a field $k$, $J\subset R$ be a homogeneous
ideal and $I\subset R$ be the Jacobian ideal of $J$, by which we always mean the ideal
$(J,I_r(\Theta))$ where $r=\hht(J)$ and $\Theta$ stands for the Jacobian matrix of a set of generators
of $J$. More precisely, if $J=(f_1,\ldots,f_s)$, then,
$$
\Theta =
\begin{bmatrix}
  \frac{\partial f_1}{\partial x_{1}} & \frac{\partial f_{2}}{\partial x_{1}} & \cdots & \frac{\partial f_{s}}{\partial x_{1}} \\
  \vdots & \vdots & & \vdots \\
 \frac{\partial f_{1}}{\partial x_{n}} & \frac{\partial f_{2}}{\partial x_{n}} & \cdots & \frac{\partial f_{s}}{\partial x_{n}} \\
\end{bmatrix}.
$$

\begin{Corollary}\label{power}
With the above assumptions and notations, if $I_r(\Theta)=\fm^r$, then the pair $J\subseteq I$ is Aluffi torsion-free.
\end{Corollary}
\demo
Let $t$ be a positive integer. Then, we have
\begin{eqnarray*}
J\cap I^t &=& J\cap (J,I_r(\Theta))^t = J\cap (J,\fm^r)^t\\
&=& J\cap (J^t, J^{t-1}\fm^r, \ldots, J\fm^{r(t-1)}) + J\cap \fm^{rt}\\
&=& J(J, \fm^r)^{t-1} + J\cap \fm^{rt} \subseteq JI^{t-1} + J\cap \fm^{rt}.
\end{eqnarray*}
The Lemma~\ref{forms} implies that $J\cap \fm^{rt}\subseteq J\fm^{r(t-1)}\subseteq JI^{t-1}$.
\qed

\bigskip

\section{Ideal of 2-minors of a $2\times n$ matrix of linear forms}
We recall the Kronecker-Weierstrass normal form of a $2\times n$ matrix of linear forms (\cite{gantmacher}).
Assume that $k$ is an algebraically closed field. Let $S$ be the polynomial
ring in variables $x_{ij}, y_{ij}, z_{ij}$ over $k$. Let $M$ be a $2\times n$ matrix of linear forms of $S$. Then, $M$ is
conjugate to a matrix obtained by concatenation of certain blocks such as
\begin{equation}\label{KW}
\left[ D_1 | \cdots | D_k | J_1 | \cdots | J_s | B_1 | \cdots | B_t \right],
\end{equation}
where $D_i$ is a ``nilpotent block" of length $n_i+ 1$:
$$
D_i = \begin{bmatrix}
x_{i1} & x_{i2} & \ldots & x_{in_i} & 0 \\
0 & x_{i1} &  \ldots & x_{i,n_i-1} & x_{in_i}
\end{bmatrix} ,
$$
$J_i$ is a ``Jordan block" of length $m_i$ with eigenvalue $\lambda_i \in k$:
$$
J_i = \begin{bmatrix}
y_{i1} & y_{i2} & \ldots & y_{im_i} \\
\lambda_i y_{i1} & y_{i1}+\lambda_i y_{i2} &\ldots &  y_{i,m_i-1} + \lambda_i y_{im_i}
\end{bmatrix} ,
$$
and $B_i$ is a ``scroll block" of length $l_i$:
$$
B_i = \begin{bmatrix}
z_{i1} & z_{i2} & \ldots & z_{i,l_i-1} & z_{il_i} \\
z_{i0} & z_{i1} &  \ldots & z_{i,l_i-2} & z_{i,l_i-1}
\end{bmatrix}.
$$
Let $I_2(M)$ be the ideal
generated by 2-minors of $M$. Since this ideal does not change under conjugation of the matrix, we will assume that $M$ is in the form of Kronecker-Weierstrass normal form.

\begin{Lemma}[\cite{rashidthesis}]\label{SJordan}
Let $M$ be a $2\times n$ matrix of linear forms in the Kronecker-Weierstrass normal form
$$
\left[ D_1 | \cdots | D_k | J_{11} | \cdots | J_{1l_1} | \cdots |J_{s_1}|\cdots | J_{sl_s} | B_1 | \cdots | B_t \right],
$$
where, each $J_{ij}$ is a Jordan block with length $p_{ij}$ and eigenvalue $\lambda_i$. Suppose that, there is at least one Jordan block with eigenvalue zero and
$$
\begin{bmatrix}
y_{1} & y_{2} & \ldots & y_{j} \\
0 & y_{1} &  \ldots &  y_{j-1}
\end{bmatrix}
$$
be the Jordan block with smallest length. Then, the ideal $(I_2(M)\colon y_1)$ is generated by all indeterminates appearing in the second row of $M$.
\end{Lemma}

\demo
Let $M'$ be the matrix obtained by deleting the column $\begin{bmatrix} y_1\\ 0 \end{bmatrix}$ and substituting $y_1$ with 0 in the matrix $M$. Denote by $J$ the ideal generated by indeterminates in the second row of $M$. Then we have the following sequence.
\begin{equation}\label{seq}
0\to \frac{S}{J}(-1) \stackrel{y_1}{\longrightarrow} \frac{S}{I_2(M)} \to \frac{S}{(I_2(M'),y_1)}\rar 0
\end{equation}
We claim that this sequence is exact. To prove it, we compare Hilbert series of them. By \cite[(2.2.3), (2.5.5)]{chun},
the Hilbert series of $S/I_2(M)$, is
$$
\frac{1}{(1-\nu)^t}\left(\frac{1+A\nu}{1-\nu} + \sum_{i=1}^s\sum_{j=1}^{l_i} \frac{p_{ij}}{(1-\nu)^{l_i-j+1}}\right) + G(\nu),
$$
where, $A=\sum_{i=1}^t m_i -1$, and $G(\nu)$ is a polynomial which is the Hilbert series of a matrix consisting of all nilpotent blocks of $M$. In the other hand,
$$
HS_{S/(I_2(M'),y_1)}(\nu) = HS_{S'/I_2(M')}(\nu),
$$
where, $S'$ is the ring $S$ without $y_1$. Since $M'$ has one column less than $M$, then
$$
HS_{S'/I_2(M')}(\nu) = \frac{1}{(1-\nu)^t}\left(\frac{1+A\nu}{1-\nu} + \sum_{i=1}^s\sum_{j=1}^{l_i} \frac{p_{ij}}{(1-\nu)^{l_i-j+1}}-\frac{\nu}{(1-\nu)^{l_1}}\right) + G(\nu).
$$
Hence,
$$
HS_{S/I_2(M)}(\nu) - HS_{S/(I_2(M'),y_1)}(\nu) = \frac{\nu}{(1-\nu)^{t+l_1}} \ ,
$$
where, $l_1$ is the number of Jordan blocks with eigenvalue zero. Note that the number of indeterminates which does not appear in the second column of $M$ is $t+l_1$. Therefore, $S/J$ is isomorphic with a polynomial ring with $t+l_1$ indeterminates. Thus, the sequence (\ref{seq}) is exact and $J=(I_2(M)\colon y_1)$.
\qed

\medskip

Note that in the above Lemma, assuming that $y_1$ is in the Jordan block with the smallest length is necessary. For example, in the matrix
$$
\begin{bmatrix}
y_1 & y_2 & w_1 & w_2 & w_3 \\
0 & y_1 & 0 & w_1 & w_2
\end{bmatrix} ,
$$
we have $y_2w_1\in I_2(M)$ but $y_2$ is not in the second row.

\begin{Proposition}\label{hight}
Let $M$ be a $2\times m$ matrix of linear forms in the Kronecker-Weierstrass normal form (\ref{KW}).
Then, hight of $I_2(M)$, the ideal generated by $2$-minors of $M$, is given in the following cases.
\begin{itemize}
\item[\rm (i)] If $M$ consists of only $k\geq 1$ nilpotent blocks, then
$$
\hht(I_2(M)) = \sum^k_{i=1} n_i.
$$
\item[\rm (ii)]
If $M$ consists of $t\geq 1$ scroll and $k\geq 0$ nilpotent blocks, then
$$
\hht(I_2(M)) = \sum^k_{i=1} n_i + \sum^t_{i=1} l_i - 1 .
$$
\item[\rm (iii)] If $M$ consists of $k\geq 0$ nilpotent, $t\geq 0$ scroll and $s\geq 1$ Jordan blocks, then
$$
\hht(I_2(M)) = \sum^k_{i=1} n_i + \sum^t_{i=1}l_i + \sum^s_{i=1}m_i  - \gamma,
$$
where, $\gamma$ is the maximum number of Jordan blocks with the same eigenvalue.
\end{itemize}
\end{Proposition}

\demo
(i) Let $M$ be of the form
$$
M = \left[\begin{array}{ccccc}
x_{1,1}& x_{1,2} & \ldots & x_{1,n_1} & 0 \\
0 & x_{1,1} &  \ldots & x_{1,n_1-1} & x_{1,n_1}
\end{array}\right| \cdots \left|
\begin{array}{ccccc}
x_{k,1} & x_{k,2} & \ldots & x_{k,n_k} & 0 \\
0 & x_{k,1} &  \ldots & x_{k,n_k-1} & x_{k,n_k}
\end{array}\right] .
$$
By \cite[P. 15]{BV},
$$
I_2(M)=\langle x_{1,1}, x_{1,2}, \ldots, x_{1,n_1}, \ldots, x_{k,1}, x_{k,2}, \ldots, x_{k,n_k}\rangle^2.
$$
Therefore, (i) is clear.

(ii) If $M$ consists of only $t$ scroll blocks. then by \cite{chun}, the Hilbert series of $S/I_2(M)$ is equal to
$$
\frac{1+(m-1)\nu}{(1-\nu)^{t+1}}.
$$
This proves the assertion in case (ii) when we have only scroll blocks.

Suppose that $M$ consists of $t\geq 1$ scroll and $k\geq 1$ nilpotent blocks. In this case, proof is by
induction on number of columns of $M$. Let $x_{11}$ be the first indeterminate in the first nilpotent block. We have the following short exact sequence:
$$
0 \to \frac{S}{I_2(M) \colon x_{11}} \stackrel{x_{11}}{\longrightarrow} \frac{S}{I_2(M)} \to \frac{S}{(I_2(M),x_{11})} \to 0.
$$
Note that,
$$
 \frac{S}{(I_2(M),x_{11})} \simeq \frac{S'}{I_2(M')},
$$
where, $M'$ is the matrix obtained by deleting the first column of $M$ and replacing 0 instead of $x_{11}$ in the second column of $M$, and $S'$ is the polynomial ring $S$ without $x_{11}$. By induction hypothesis, there is
$h'(\nu)\in {\mathbb Z}[\nu]$ such that, the Hilbert series of $S'/I_2(M')$ is of the form
$$
\frac{h'(\nu)}{(1-\nu)^{c-1-(\delta-2)}} = \frac{h'(\nu)}{(1-\nu)^{c-(\delta-1)}},
$$
where, $\delta=\sum^k_{i=1} n_i + \sum^t_{i=1} l_i$, and $c$ is number of all indeterminates appearing in $M$.

If $l_i\geq 3$, for $i=1,\ldots, t$, then, the ideal $I_2(M)\colon x_{11}$ is generated by all indeterminates.
If for some $1\leq i\leq t$, $1\leq l_i\leq 2$, then, $z_{i,l_i}^u\in I_2(M)\colon x_{11}$, for some positive
integer $u$.
Since the ideal $(I_2(M)\colon x_{11})$ is zero dimensional, therefore, the Hilbert series of
$S/(I_2(M)\colon x_{11})(-1)$ is simply $\nu h(\nu)$ for some $h(\nu)\in {\mathbb Z}[\nu]$.
By using the above short exact sequence and additive property of Hilbert series, we
obtain the Hilbert series of $S/I_2(M)$:
$$
H_{S/I_2(M)}(\nu) = \frac{\nu h(\nu)(1-\nu)^{c-(\delta-1)} + h'(\nu)}{(1-\nu)^{c-(\delta - 1)}} .
$$
In this fraction, the numerator is not divisible by $(1-\nu)$. Therefore, dimension of $S/I_2(M)$ is $c-(\delta-1)$
and hight of $I_2(M)$ is $\delta-1$. This completes the proof of case (ii).

(iii) Suppose that $M$ has $s\geq 1$ Jordan blocks. Also in this case, the proof is by
induction on number of columns of $M$.
Let $\gamma$ be the maximum number of Jordan blocks with the same eigenvalues
$\lambda$. After some suitable elementary column and row operations, we obtain a matrix conjugate to $M$ such that
lengths and types of all blocks are preserved and the blocks with eigenvalue $\lambda$ have become to blocks
with eigenvalue zero (for details, see the proof of the main theorem in \cite{CJ}).
Let $y_{11}$ be the first
indeterminate in the smallest Jordan block with eigenvalue zero. The above short exact sequence is valid if we
substitute $x_{11}$
by $y_{11}$. In this case,                                                                                                   $$
 \frac{S}{(I_2(M),y_{11})} \simeq \frac{S'}{I_2(M')},
$$
where, $M'$ is the matrix obtained by  $M$ deleting first column and replacing 0 instead of $y_{11}$, and $S'$ is the polynomial ring $S$
without $y_{11}$. By induction hypothesis, there is
$h'(\nu)\in {\mathbb Z}[\nu]$ such that, the Hilbert series of $S'/I_2(M')$ is of the form
$$
\frac{h'(\nu)}{(1-\nu)^{c-1-(\delta-\gamma-1)}} = \frac{h'(\nu)}{(1-\nu)^{c-(\delta-\gamma)}},                                               $$
where, $\delta=\sum^k_{i=1} n_i + \sum^t_{i=1} l_i + \sum_{i=1}^s m_i$. Note that, the ideal $(I_2(M)\colon y_{11})$ is generated by all indeterminates appearing in the second row of $M$ (\cite{ZZ}).  The number of indeterminates appearing
in the second row is $\delta-\gamma$. Therefore, the Hilbert series of
 $S/(I_2(M)\colon y_{11})(-1)$ is
$$
\frac{\nu}{(1-\nu)^{c-(\delta - \gamma)}}.
$$
The Hilbert series of $S/I_2(M)$ is
$$
H_{S/I_2(M)}(\nu) = \frac{\nu + h'(\nu)}{(1-\nu)^{c-(\delta - \gamma)}} .
$$
Therefore, dimension of $S/I_2(M)$ is $c-(\delta-\gamma)$
and hight of $I_2(M)$ is $\delta-\gamma$.
\qed

\bigskip

\begin{Theorem}\label{matrix}
Let $M$ be a $2\times n$ matrix of linear forms in a polynomial ring $S$ over an algebraically closed field $k$. Suppose that $I_2(M)$ has codimension  $r>1$. Denote by $\Theta$ the Jacobian matrix of $I_2(M)$.  Then, the following conditions are equivalent.
\begin{itemize}
\item[\rm (a)]
$I_r(\Theta) = \fm^r$.
\item[\rm (b)]
The Kronecker-Weierstrass normal form of $M$ does not have any Jordan block, or it consists of only some nilpotent blocks and some Jordan blocks of length 1.
\item[\rm (c)]
The pair $I_2(M)\subseteq (I_2(M), I_r(\Theta))$ is Aluffi torsion-free.
\end{itemize}
Where $\fm$ is the irrelevant maximal ideal of $S$ and $I_r(\Theta)$ is the ideal generated by $r$-minors of $\Theta$.

\end{Theorem}

\demo
(a $\Rightarrow$ b) Let $M$ be a matrix which has at least one Jordan block.
Suppose that $\gamma$ is the maximum number of Jordan blocks with the same eigenvalue
$\lambda$. As stated in the proof of Proposition~\ref{hight}, we may assume that $\lambda$ is zero.
Let the block $J_1$ be one of the Jordan blocks with length greater than 1 and eigenvalue zero. It is in the form:
$$
\begin{bmatrix}
y_{1,1} & y_{1,2} & \ldots & y_{1,m_1} \\
0 & y_{1,1} &\ldots &  y_{1,m_1-1}
\end{bmatrix} .
$$
By the Proposition~\ref{hight}, height of $I_2(M)$ is
$$
r=\sum_1^k
 n_i + \sum_1^s m_i + \sum_1^t l_i -\gamma.
$$
But, the variable $y_{1,m_1}$ appears only in $r-1$ quadratic forms in generators of $I_2(M)$ and therefore, it appears only in $r-1$ rows in the Jacobian matrix of $I_2(M)$. This is enough to know that $y_{1,m_1}^r$ is not in the ideal of $r$-minors of the Jacobian matrix of $I_2(M)$ and then, $I_r(\Theta) \neq \fm^r$.

If there is no any Jordan block of length greater than 1, and there is at least one scroll block and some at least one Jordan block of length 1, then the variable $z_{1l_1}$ appears in $r-1$ quadratic forms in $I_2(M)$ and the same argument as above shows that $z_{1l_1}^r\not\in I_r(\Theta)$.

\medskip

(b $\Rightarrow$ a)
If the Kronecker-Weierstrass normal form of $M$ does not have any Jordan block, then, $M$ falls within one of the following cases.
\begin{itemize}
\item[\rm (i)] $M$ has only scroll blocks,
\item[\rm (ii)] $M$ has only nilpotent  blocks,
\item[\rm (iii)] $M$ has nilpotent and scroll blocks.
\item[\rm (iv)] $M$ has nilpotent and Jordan blocks of length 1.
\end{itemize}
In each case, we show that $I_r(\Theta) = \fm^r$.

\medskip

\underline{{\it Case i}}.
First assume that there is only one scroll block:
$$
M = \left[\begin{array}{ccccc}
z_{1} & z_{2} & \ldots & z_{m-1} & z_{m} \\
z_{0} & z_{1} &  \ldots & z_{m-2} & z_{m-1}
\end{array}\right] .
$$
We prove that for each monomial of degree $m-1$, there is a $(m-1)$-minor of $\Theta$ such that the monomial is initial term of the minor with lexicographic order.
By
$$A=[z_{i_1} \  z_{i_2} \ldots z_{i_r} \ | \ (c_{11},c_{12}) \ (c_{21},c_{22}) \ldots (c_{r1},c_{r2})]$$
we mean the $r$-minor of $\Theta$ such that the entry $[A]_{kl}= \partial f_{(c_{l1},c_{l2})}/\partial z_k$, where, $f_{(c_{l1},c_{l2})}$ is the 2-minor of $M$ obtained by columns $c_{l1}$ and $c_{l2}$. The following equations are clear.
\begin{eqnarray*}
z_0^{m-1} &=&  [z_{2} \ z_{3} \ldots z_{m} \ | \ (1,2) \ (1,3) \ldots (1,m)]\\
z_1^{m-1} &=&  [z_{1} \ z_{3} \ldots z_{m} \ | \ (1,2) \ (2,3) \ldots (2,m)]\\
 & \vdots &\\
z_i^{m-1} &=&  [z_{i} \ z_{i-2} \ldots z_{0} \ z_{i+2} \ldots z_{m}| (i,i+1) \ (i-1,i) \ldots (1,i) \ (i+1,i+2) \ldots (i+1,m)]\\
 & \vdots &\\
z_{m-1}^{m-1} &=& [z_{m-1} \ z_{m-3} \ldots z_{0} \ | \ (m-1,m) \ (m-2,m-1) \ldots  (1,m-1) ]\\
z_{m}^{m-1} &=& [z_{m-2} \ z_{m-3} \ldots z_{0} \ | \ (m-1,m) \ (m-2,m) \ldots (1,m) ]
\end{eqnarray*}
All above minors are upper triangular.

Let $z_{j_1}^{a_1}z_{j_2}^{a_2}\cdots z_{j_s}^{a_s}$ be a given monomial of degree $r=m-1$. Take the minor
$$
[z_{d_1} \  \ldots z_{d_r} \ | \ (h_{11},h_{12}) \ (h_{21},h_{22}) \ldots (h_{r1},h_{r2})]
$$
such that $a_1$ of $z_i$'s are first $a_1$ entries of the minor of $z_{j_1}^{m-1}$. Then, for the succeeding $a_2$ of $z_i$'s choose first $a_2$ entries of the minor of $z_{j_2}^{m-1}$, which they are not appeared in the previous chooses and also the columns are not repeated. Continuing this process, we get a minor which its main diagonal is the given monomial $z_{j_1}^{a_1}z_{j_2}^{a_2}\cdots z_{j_s}^{a_s}$ and this monomial is initial of the minor. To show the last statement, note that entries below the main diagonal do not effect the initialness of the main diagonal. Example~\ref{examscroll} illustrates concretely this argument.

Now let $M$ be of the form
$$
M = \left[\begin{array}{ccccc}
z_{1,1} & z_{1,2} & \ldots & z_{1,l_1-1} & z_{1,l_1} \\
z_{1,0} & z_{1,1} &  \ldots & z_{1,l_1-2} & z_{1,l_1-1}
\end{array}\right| \cdots \left|
\begin{array}{ccccc}
z_{c,1} & z_{c,2} & \ldots & z_{c,l_c-1} & z_{c,l_c} \\
z_{c,0} & z_{c,1} &  \ldots & z_{c,l_c-2} & z_{c,l_c-1}
\end{array}\right] .
$$
First consider the lexicographic order on terms of $S$ with respect to $z_{1,0}>z_{1,1}>\cdots>z_{c,l_c}$ and write the generators of $I_2(M)$ with this order:
$$
I_2(M)=(f_1,\ldots,f_t,f_{t+1},\ldots,f_k),
$$
where, $z_{1,0}$ appears in $f_1, \ldots, f_t$ and does not appear in $f_{t+1},\ldots,f_k$. Then, the Jacobian matrix of $I_2(M)$ is of the form:
\begin{equation}\label{scroll}
\Theta=\left[
\begin{array}{cccccccc|ccccc}
-z_{1,2} & -z_{1,3} & \cdots & -z_{1,l_1} & \cdots &  -z_{c,1} & \cdots & -z_{c,l_c} & 0 & 0 & \cdots & 0 & 0 \\
\hline
&&&&&&&&&&&& \\
&&&&\ast &&&&&& \Theta' && \\
&&&&&&&&&&&&
\end{array}
\right]
\end{equation}
In this matrix, the block $\Theta'$ is Jacobian matrix of $I_2(M')$ where $M'$ is a matrix obtained by deleting first column of $M$. By induction on number of columns of $M$, we have $I_{r-1}(\Theta') = {\fm'}^{r-1}$, where $\fm'$ is the ideal $\fm$ without $z_{1,0}$. By the form of $\Theta$, it is clear that
$$
z_{i,j}\langle z_{1,1}, z_{1,2}, \ldots, z_{c,l_c}\rangle^{r-1}\subseteq I_r(\Theta), \ \ \ 1\leq i\leq c, \ \ \ 1\leq j\leq l_i, \ \  \ (i,j)\neq (1,1).
$$
Therefore, $$\langle z_{1,2}, \ldots, z_{1,l_1}, \ldots, z_{c,1}, \ldots, z_{c,l_c-1}, z_{c,l_c}\rangle^r\subseteq I_r(\Theta). $$
In other hand, if we assume the degree reverse lexicographic order with respect to $z_{1,0}>z_{1,1}>\cdots>z_{c,l_c}$, the Jacobian matrix of $I_2(M)$ is of the form:
$$
\left[
\begin{array}{ccccc|cccccccc}
&&&&&&&&&&&& \\
&& \Theta'' &&&&&&& \ast &&& \\
&&&&&&&&&&&& \\
\hline
0 & 0 & \cdots & 0 & 0 & -z_{1,0} & -z_{1,1} & \cdots & -z_{1,l_1-1} & \cdots &  -z_{c,0} & \cdots & -z_{c,l_c-2}
\end{array}
\right].
$$
Where, the block $\Theta''$ is Jacobian matrix of $I_2(M'')$ where $M''$ is a matrix obtained by deleting the last column of $M$. Note that the latter matrix is obtained by some changes of columns of the matrix $\Theta$. Again  by induction on number of columns of the matrix $M$, we have $I_{r-1}(\Theta'') = {\fm''}^{r-1}$, where $\fm''$ is the ideal $\fm$ without $z_{c,l_c}$. Then, it is clear that
$$
z_{i,j}\langle z_{1,0}, z_{1,2}, \ldots, z_{c,l_c-1}\rangle^{r-1}\subseteq I_r(\Theta), \ \ \ 1\leq i\leq c,\ \ \  0\leq j\leq l_i-1,\ \  \ (i,j)\neq (c,l_c-1).
$$
Therefore, $$\langle z_{1,0}, \ldots, z_{1,l_1-1}, \ldots, z_{c,0}, \ldots, z_{c,l_c-1}, z_{c,l_c-2}\rangle^r\subseteq I_r(\Theta). $$
Changing the first and last blocks of $M$ and repeating the above argument, completes the proof in this case.

\medskip

\underline{{\it Case ii}}.

If the matrix $M$ consists of only nilpotent blocks, then by Proposition~\ref{hight},  $I_2(M)=\fm^2$ and  clearly $I_r(\Theta) = \fm^r$.

\medskip

\underline{{\it Case iii}}. Let $M$ be a matrix obtained by concatenation of some scroll blocks and some nilpotent blocks:
$$
M= [ D_1 | \cdots | D_r | B_1 | \cdots | B_t ].
$$
Let $x_{11}$ be the first entry of the first nilpotent block $D_1$. Then, $x_{11}^2\in I_2(M)$ and with the same method of case (i), the Jacobian matrix of $I_2(M)$ will be in the form of (\ref{scroll}) with all indeterminates  appearing in the top-left block. Using induction on number of columns of $M$ proves the theorem in this case.

\medskip

\underline{{\it Case iv}}. Let $M$ be a matrix consisting of $k\geq 0$ nilpotent blocks and $s$ Jordan blocks:
$$
M = \left[\begin{array}{ccc|c}
x_{1,1}& \ldots &  0 & \ldots\\
0 &   \ldots  & x_{1,n_1} & \ldots
\end{array}
\begin{array}{|ccc}
x_{k,1} &  \ldots  & 0 \\
0 &  \ldots &  x_{k,n_k}
\end{array}
\begin{array}{|c|c|c|c|c|c}
y_{1} &  \ldots  & y_{\gamma} & y_{\gamma+1} & \ldots & y_s\\
0 &  \ldots &  0 & \lambda_1 y_{\gamma+1} &  \ldots & \lambda_s y_s
\end{array}
\right] .
$$
If $k>0$, then the same argument as case (iii) concludes case (iv). If there is no any nilpotent block, take $y_1$,  $y_2$, $y_s$ and use the induction argument as in case (iii).

\medskip

(a $\Rightarrow$ c) It follows from Corollary~\ref{power}.

(c $\Rightarrow$ b) Let $M$ has Jordan blocks of length greater than 1. In this case,  it is clear that $f=(y_{1,1}y_{1,m_1}^{r-1})\in I_r(\Theta)\setminus I_2(M)$ but,
$f^2\in I_2(M)$. Therefore, $f^2\in I_2(M)\cap I_r(\Theta)^2$ but, $f^2\not\in I_2(M)I_r(\Theta)$.
Let $M$ has $t>0$ scroll blocks and $s>0$ Jordan blocks of length 1. Then, $f_1=z_{1l_1-1}z_{1l_1}^{r-1}$ and $f_2=z_{1l_1}^{r-1}y_1$ are in $I_r(\Theta)$, but they are not in $I_2(M)$. In other hand, $f_1f_2\in I_2(M)\cap I_r(\Theta)^2$ but, $f_1f_2\not\in I_2(M)I_r(\Theta)$.
\qed

\begin{Example}\label{examscroll}\rm
Let $M$ be the matrix
$$
\left[\begin{array}{ccccccc}
z_{1} & z_{2} & z_3 & z_4 & z_5 & z_6 & z_7\\
z_{0} & z_{1} & z_2 & z_3 & z_4 & z_5 & z_6
\end{array}\right] .
$$
Following is illustration of some monomials as initials of minors.
\begin{eqnarray*}
z_0^{6} &=&  [z_{2} \ z_{3} \ z_4 \ z_5 \ z_6 \  z_{7} \ | \ (1,2) \ (1,3) \ (1,4) \ (1,5) \ (1,6) \ (1,7)]\\
z_1^{6} &=&  [z_{1} \ z_{3} \ z_4 \ z_5 \ z_6 \  z_{7} \ | \ (1,2) \ (2,3) \ (2,4) \ (2,5) \ (2,6) \ (2,7)]\\
z_2^{6} &=&  [z_{2} \ z_{0} \ z_4 \ z_5 \ z_6 \  z_{7} \ | \ (2,3) \ (1,2) \ (3,4) \ (3,5) \ (3,6) \ (3,7)]\\
z_3^{6} &=&  [z_{3} \ z_{1} \ z_0 \ z_5 \ z_6 \  z_{7} \ | \ (3,4) \ (2,3) \ (1,3) \ (4,5) \ (4,6) \ (4,7)]\\
z_4^{6} &=&  [z_{4} \ z_{2} \ z_1 \ z_0 \ z_6 \  z_{7} \ | \ (4,5) \ (3,4) \ (2,4) \ (1,4) \ (5,6) \ (5,7)]\\
z_5^{6} &=&  [z_{5} \ z_{3} \ z_2 \ z_1 \ z_0 \  z_{7} \ | \ (5,6) \ (4,5) \ (3,5) \ (2,5) \ (1,6) \ (6,7)]\\
z_6^{6} &=&  [z_{6} \ z_{4} \ z_3 \ z_2 \ z_1 \  z_{0} \ | \ (6,7) \ (5,6) \ (4,6) \ (3,6) \ (2,6) \ (1,6)]\\
z_7^{6} &=&  [z_{5} \ z_{4} \ z_3 \ z_2 \ z_1 \  z_{0} \ | \ (6,7) \ (5,7) \ (4,7) \ (3,7) \ (2,7) \ (1,7)]
\end{eqnarray*}
$$
z_0z_1^2z_4z_7^2 = \mbox{In}([z_{2} \ z_{3} \ z_4 \ z_0 \ z_5 \  z_{1} \ | \ (1,2) \ (2,3) \ (2,4) \ (1,4) \ (6,7) \ (2,7)]).
$$
The sub-matrix corresponding to the last monomial is:
$$
\begin{bmatrix}
-z_0 & 2z_2 & z_3 & 0 & 0 & z_6 \\
0 & -z_1 & z_2 & z_1 & 0 & 0 \\
0 & 0 & -z_1 & -z_0 & 0 & 0 \\
-z_2 & 0 & 0 & -z_4 & 0 & 0 \\
0 & 0 & 0 & 0 & -z_7 & 0 \\
2z_1 & -z_3 & z_4 & z_3 & 0 & -z_7
\end{bmatrix} .
$$
\end{Example}

\begin{Remark}\rm
Let $X\subset \mathbb{P}_k^{n}$ be a projective algebraic set of dimension $d$ with defining ideal $I_2(M)$ where $M$ is a matrix of linear forms in
$k[x_0,\ldots,x_n]$ and $k$ is algebraically closed. Theorem~\ref{matrix} gives a criterion to check nonsingularity of $X$, that is, the Kronecker-Weierstrass normal form of $M$ does not have any Jordan block, or it consists of only some nilpotent blocks and some Jordan blocks of length 1 if and only if $X$ is non-singular.

Note that by proof of the above theorem, in case that $M$ does not have any Jordan block, then the ideal $I_r(\Theta)$ is $\fm$-primary but, in the case that $M$ has Jordan blocks, it is not $\fm$-primary. This means that the following conditions are equivalent.
\begin{itemize}
\item[\rm (a)]
$I_r(\Theta) = \fm^r$.
\item[\rm (b)]
$I_r(\Theta)$ is $\fm$-primary.
\end{itemize}
\end{Remark}

This Remark initiates the following conjecture.

\begin{Conjecture}\label{conj1}\rm
Let $J$ denote the ideal generated by quadrics in a polynomial ring $S$, such that $r=\hht(J)\geq 2$. Then, the following conditions are equivalent.
\begin{itemize}
\item[\rm (a)]
$I_r(\Theta) = \fm^r$.
\item[\rm (b)]
$I_r(\Theta)$ is $\fm$-primary.
\end{itemize}
Where, $\Theta$ is the Jacobian matrix of $J$ and $\fm$ is the irrelevant maximal ideal of $S$.
\end{Conjecture}

\begin{Corollary}\label{HB}
Let $J\subset R=k[x_1,\ldots,x_n]$ denote a codimension $2$ ideal generated by $3$ quadrics with the following free resolution:
$$
0\to R^2 \to R^3 \to J \to 0.
$$
Let $I_2(\Theta)$ denote the ideal generated by the $2$-minors of the Jacobian matrix $\Theta$ of the generators of
$J$. If $I_2(\Theta)$ is $\fm=(x_1,\ldots,x_n)$-primary, then the pair $J\subset (J,I_2(\Theta))$ is Aluffi torsion-free. In particular $V(J)\subseteq {\mathbb P}^{n-1}$ is nonsingular.
\end{Corollary}
\demo
By the Hilbert-Burch theorem, $J$ is generated by $2$-minors of the syzygy matrix $M$ of $J$. By assumption, the transpose of $M$ is a $2\times 3$ matrix of linear forms in $R$. Since $I_2(\Theta)$ is $\fm$-primary, Theorem ~\ref{matrix} implies that the Kronecker-Weierstrass normal form of $M$ does not have Jordan block and $I_2(\Theta)=\fm^2$. Then by Corollary~\ref{power}, the pair $J\subset (J,I_2(\Theta))$ is Aluffi torsion-free. Since the Jacobian ideal has codimension $n$, then the additional assertion at the end of the statement is clear.
\qed

\medskip

Recall that a $n\times n$ (generic) Hankel matrix is of the form
$$
H=
\begin{bmatrix}
  x_1 & x_2 &  \ldots & x_{n-1} & x_n\\
  x_2 & x_3 &  \ldots & x_n & x_{n+1}\\
   \vdots & \vdots &  & \vdots & \vdots \\
    x_{n-1} & x_n & \ldots & x_{2n-3} & x_{2n-2}\\
  x_n & x_{n+1} & \ldots & x_{2n-2} & x_{2n-1}
\end{bmatrix} ,
$$
and a generalized Hankel matrix is concatenation of some Hankel matrices (with different indeterminates).

\begin{Corollary}
Let $J$ be the ideal of $2$-minors of
 a generalized Hankel matrix. Then, the pair $J\subseteq (J, I_r(\Theta))$ is Aluffi torsion-free.
\end{Corollary}

\demo
By \cite[Theorem 2.2]{ZU}, $J$ is generated by 2-minors of a $2\times m$ matrix which has only scroll blocks. Now, use the  Theorem~\ref{matrix} to complete the proof.
\qed

\begin{Examples}\label{rational}\rm
\begin{itemize}
\item[\rm (i)] The rational normal scroll in ${\mathbb P}^d_k$, could be realized as the variety of the ideal $J$ generated by 2-minors of a matrix consisting only scroll blocks \cite{harris}.
If $I$ is the Jacobian ideal of $J$, then by Theorem~\ref{matrix}, the pair $J\subseteq I$ is Aluffi torsion-free.
\item[\rm (ii)] Consider the rational map $F : {\mathbb P}^2_k\dashrightarrow {\mathbb P}^4_k$ given by
$$
F(y_0:y_1:y_2) = (y_0^2 : y_1^2 : y_0y_1 : y_0y_2 : y_1y_2).
$$
The image of this map is given by the ideal
$$
J = \langle x_2^2 - x_0x_1, x_2x_3 - x_0x_4, x_2x_4 - x_1x_3\rangle .
$$
Note that $J$ is generated by 2-minors of the matrix
$$
\begin{bmatrix}
  x_2 & x_1 & x_4 \\
  x_0 & x_2 & x_3
\end{bmatrix},
$$
which consists of two scroll blocks. Therefore, the pair $J\subseteq I$ is Aluffi torsion-free.
\end{itemize}
\end{Examples}

\bigskip

\section{Edge ideal of a graph}
Let $I$ be a monomial ideal in the polynomial ring $k[x_1,\ldots,x_n]$. It is known that the ideal of $r$-minors of the Jacobian matrix of $I$ is again a monomial ideal (see \cite{Simis} and \cite{Thesis2}). We provide another simple proof for this fact in Lemma~\ref{jacobian}.

Let $M$ be a $m\times n$ matrix and $1\leq r\leq \min\{m,n\}$ be an integer. A transversal of length $r$ in $M$ or an $r$-transversal of $M$ is a product of $r$ entries of $M$ with different rows and columns. In other words, an $r$-transversal of $M$ is product of entries of
 the main diagonal of an $r\times r$ sub-matrix of $M$ after suitable changes of columns and rows.

\begin{Lemma}\label{jacobian}
Let $I$ be an ideal of $k[x_1,\ldots,x_n]$ generated by monomials $m_1, \ldots, m_s$. Let $\Theta$ be the Jacobian matrix of $I$ and $1\leq r\leq \min\{n,s\}$. Then, any $r$-minor of $\Theta$ is a monomial.
\end{Lemma}
\demo
Let $f=[a_1,\ldots,a_r|b_1,\ldots,b_r]$ represent an $r$-minor of $\Theta$. That is, $1\leq a_1<a_2<\cdots<a_r\leq n$ are rows and  $1\leq b_1<b_2<\cdots<b_r\leq s$ are columns of the matrix $\Theta$ appearing in the chosen $r$-minor.
The corresponding sub-matrix is:
$$
\begin{bmatrix}
  \frac{\partial m_{a_1}}{\partial x_{b_1}} & \frac{\partial m_{a_2}}{\partial x_{b_1}} & \cdots & \frac{\partial m_{a_r}}{\partial x_{b_1}} \\
  \vdots & \vdots & & \vdots \\
 \frac{\partial m_{a_1}}{\partial x_{b_r}} & \frac{\partial m_{a_2}}{\partial x_{b_r}} & \cdots & \frac{\partial m_{a_r}}{\partial x_{b_r}} \\
\end{bmatrix}
$$

Note that, any term of $f$ is an $r$-transversal. This term is nonzero if in any factor $\frac{\partial m_{a_i}}{\partial x_{b_j}}$ of it, $m_{a_i}$ is divisible by $x_{b_j}$ and in this case,  $\frac{\partial m_{a_i}}{\partial x_{b_j}}= \gamma \frac{m_{a_i}}{x_{b_j}}$, where the integer $\gamma$ is the highest power of $x_{b_j}$ appearing in $m_{a_i}$. Therefore, any nonzero term of $f$ is of the form:
$$
\beta \frac{m_{a_1} \cdots m_{a_r}}{x_{b_1}\cdots x_{b_r}},
$$
where $\beta$ is an integer. The minor $f$ is sum of the same monomials with possibly different coefficients and therefore, it is a monomial.
\qed

\medskip

Recall that for a finite simple graph $G$ with vertex set $V(G)=\{v_1,\ldots,v_n\}$, an ideal $I(G)$ in the ring $k[x_1,\ldots,x_n]$ is corresponded which is generated by all square-free quadratic monomials $x_ix_j$ provided that $\{v_i,v_j\}$ is an edge in $G$. This ideal is called the edge ideal of $G$. Let $v$ be a vertex in $G$. Degree of $v$ is number of all vertices adjacent to $v$. For a subset $A$ of $V(G)$, the set of all vertices adjacent to some vertices in $A$ is called neighborhood of $A$ and denoted by $N(A)$. A subset $B$ of vertices of $G$ is called an independent set if there is no any edge between each two vertices of $B$.
A matching in $G$ is a subset of edges of $G$ such that there is no any common vertex between any two of them. In this section, we identify any edge $v_i$ with the corresponding indeterminate $x_i$.

\begin{Lemma}\label{maxminor}
Let $G$ be a graph with $n$ vertices, $I(G)$ edge ideal of $G$ and $\Theta$ the Jacobian matrix of $I(G)$. Let $g\in k[x_1,\ldots,x_n]$ be a monomial  and  $r$ a positive  integer. The following conditions are equivalent.
\begin{itemize}
\item[\rm (i)] $g$ is a $r$-transversal of $\Theta$.
\item[\rm (ii)] There are $r$ different edges $e_1=\{x_{1_1}, x_{1_2}\}, \ldots, e_r=\{x_{r_1},x_{r_2}\}$ such that  vertices $x_{1_1}, \ldots, x_{r_1}$ are different and $g= x_{1_{2}} \cdots x_{r_{2}}$.
 \end{itemize}
Moreover, let  the set $\{x_{i_1}, \ldots, x_{i_s}\}$ is independent. Then there is a $r$-transversal of the form  $g=x_{i_1}^{\alpha_1}\cdots x_{i_s}^{\alpha_s}$ with $0\leq \alpha_j\leq \deg(x_{i_j})$ for $1\leq j\leq s$ and $\sum \alpha_j=r$, if and only if $|N(\{x_{i_1},  \ldots, x_{i_s}\})|\geq r$.
\end{Lemma}
\demo
Generators of the ideal $I(G)$  are of the form $x_ix_j$ where $\{x_i,x_j\}$ is an edge in $G$ and each entry of the Jacobian matrix $\Theta$ is zero or of the form $x_i$ where $x_i$ is belonging to an edge in $G$. Equivalence of (i) and (ii) is clear by definition of $r$-transversal of $\Theta$.

By the Lemma~\ref{jacobian}, any $r$-transversal of $\Theta$ is a monomial of degree $r$. Let $g=x_{i_1}^{\alpha_1} \cdots x_{i_s}^{\alpha_s}$ be a $r$-transversal of $\Theta$. It means that there is a $r\times r$ sub-matrix of $\Theta$, which admits $b_1$ times $x_{i_1}$, \dots, and $b_s$ times $x_{i_s}$ in different rows and columns. In the matrix $\Theta$, the entry $x_{i_j}$ appears exactly $\deg(v_{i_j})$ times. Therefore $\alpha_j\leq \deg(x_{i_j})$ for each $1\leq j\leq s$. Moreover, if $A = \{x_{i_1}, \ldots, x_{i_s}\}$ is an independent set of vertices, then the set $N(A)$ contains vertices which are adjacent to some vertices in $A$ and there are $|N(A)|$ different edges between $A$ and $N(A)$ with different ends in $B$. Now, it is clear that there is a $r$-transversal of the form $g=x_{i_1}^{\alpha_1}\cdots x_{i_s}^{\alpha_s}$ with $0\leq \alpha_j\leq \deg(x_{i_j})$ for $1\leq j\leq s$ and $\sum \alpha_j=r$, if and only if $|N(A)|\geq r$.
\qed

\medskip

We say that a graph $G$ is Aluffi torsion-free if the pair $I(G)\subseteq(I(G),I_r(\Theta))$ is Aluffi torsion-free, where $r$ is height of $I(G)$ and $\Theta$ is Jacobian matrix of $I(G)$.

\begin{Theorem}\label{graph}
Let  $G$ be a graph and $\hht(I(G))=r>1$. Then $G$ is not Aluffi torsion-free if and only if there are adjacent vertices $x_1, x_2$ and other vertices $x_{i_1}, \ldots, x_{i_s}$ for some integer $s\geq 1$, such that
\begin{itemize}
\item[\rm (i)] The sets $\{x_1, x_{i_1}, \ldots, x_{i_s}\}$ and $\{x_2, x_{i_1}, \ldots, x_{i_s}\}$ both are independent, and
\item[\rm (ii)] $|N(\{x_{i_1}, \ldots, x_{i_s}\})|= r-1$.
\end{itemize}
\end{Theorem}
\demo
Let $G$ be not Aluffi torsion-free. Then, there is an integer $t\geq 2$ such that
\begin{equation}\label{equality}
I(G)\cap (I(G),I_r(\Theta))^t\neq I(G)(I(G),I_r(\Theta))^{t-1}.
\end{equation}
Note that the right hand side is always a subset of the left hand side and it is enough to check the reverse inclusion. Let $g$ be a monomial in left hand side which is not in right hand side of (\ref{equality}). Then $g=g_1 \cdots g_t$ such that $g_i\in (I(G),I_r(\Theta))$. If for some $1\leq i\leq t$, $g_i\in I(G)$, then $g=g_i(g_1\cdots g_{i-1}g_{i+1}\cdots g_t)\in I(G)(I(G),I_r(\Theta))^{t-1}$, which is a contradiction.
Note that a $r$-transversal $g_i$ belongs to $I(G)$ if and only if the set of vertices appearing in $g_i$ is not independent.

The monomial $g$ is in $I(G)$ then there are adjacent vertices $x_k, x_l$ such that $x_kx_l| g$, but $x_kx_l\nmid g_i$ for each $i=1,\ldots,t$. Without loss of generality, let $x_k|g_1$ and $x_l|g_2$. In this situation, $g_1g_2\in I(G)\cap (I(G),I_r(\Theta))^2$. If $g_1g_2\in I(G)(I(G),I_r(\Theta))$, then $g_3g_4\cdots g_t\in (I(G),I_r(\Theta))^{t-2}$ and $g\in I(G)(I(G),I_r(\Theta))^{t-1}$ which is again a contradiction. Therefore, we may assume that $g=g_1g_2\in I(G)\cap (I(G),I_r(\Theta))^2\setminus I(G)(I(G),I_r(\Theta))$ and $g_i\in I_r(\Theta)\setminus I(G)$ for $i=1,2$. Moreover,
$x_1|g_1$, $x_2|g_2$ and $x_1$ is adjacent to $x_2$.

Assume that $g_1= x_1 x_{i_1}^{\alpha_1}\cdots x_{i_s}^{\alpha_s}$ and $g_2= x_2 x_{j_1}^{\beta_1}\cdots x_{j_t}^{\beta_t}$, such that $\sum \alpha_i = \sum \beta_j = r-1$ and both sets $A=\{x_1, x_{i_1}, \ldots, x_{i_s}\}$ and $B=\{x_2, x_{j_1}, \ldots, x_{j_t}\}$ are independent.
If the set $\{x_{i_1}, \ldots, x_{i_s}, x_{j_1}, \ldots, x_{j_t}\}$ is dependent, then $g_1g_2 \in (I(G))^2\subseteq I(G)(I(G),I_r(\Theta))$, a contradiction. By the same argument, it is not possible that $x_1$ is adjacent to some vertex in $B\setminus\{x_2\}$ and simultaneously $x_2$ is adjacent to some vertex in $A\setminus\{x_1\}$. Assume that $x_2$ is not adjacent to any vertex in $A\setminus\{x_1\}$. We claim that the vertices $x_1$, $x_2$ and $x_{i_1}, \ldots, x_{i_s}$ satisfy conditions (i) and (ii).

By the procedure of the above argument, the vertices  $x_1$, $x_2$ and $x_{i_1}, \ldots, x_{i_s}$ clearly satisfy conditions (i). In other hand, $x_{i_1}^{\alpha_1}\cdots x_{i_s}^{\alpha_s}$ is a $(r-1)$-transversal of $\Theta$ and by Lemma~\ref{maxminor}, $|N(\{x_{i_1}, \ldots, x_{i_s}\})|\geq r-1$. We know that $x_{i_1}^{\alpha_1}\cdots x_{i_s}^{\alpha_s} x_{j_1}^{\beta_1}\cdots x_{j_t}^{\beta_t}$ is not in $I_r(\Theta)$ and thus there is no any $r$-transversal of $\Theta$ dividing it. This means that for any subset $C$ of $\{x_{i_1}, \ldots, x_{i_s}, x_{j_1}, \ldots,  x_{j_t}\}$, $|N(C)|<r$. Therefore $|N(\{x_{i_1}, \ldots, x_{i_s}\})| = r-1$, as required.

Conversely, let there are vertices $x_1$, $x_2$ and $x_{i_1}, \ldots, x_{i_s}$ satisfying conditions (i) and (ii). Let $g_3=x_2x_{i_1}\cdots x_{i_s}$. Then $g_3$ is a $r$-transversal of $\Theta$ and $g_1g_3\in I(G)\cap (I(G),I_r(\Theta))^2$. By Lemma~\ref{maxminor}, condition (ii) guarantees that $g_1g_3/x_1x_2\not\in (I(G),I_r(\Theta))$. Therefore, $G$ is not Aluffi torsion-free.
\qed

\begin{Examples}\label{ex3.8}\rm
\begin{itemize}
\item[\rm (i)]
A complete graph $K_n$ for $n>2$ is Aluffi torsion-free. Because all vertices are adjacent to each other and there is no any vertex satisfying condition (i) of the above theorem.
\item[\rm (ii)] A complete $r$-partite graph is Aluffi torsion-free. In contrary if it is not Aluffi torsion-free, then, there are two adjacent vertices $v_1, v_2$ and at least one another vertex $w$ which is adjacent to none of $v_1$ and $v_2$. In this case, $v_1$ and $w$ belongs to the same part and also $v_2$ and $w$ belongs to the same part. Therefore $v_1$ and $v_2$ are in the same part which is a contradiction.
\item[\rm (iii)] A complete graph minus edges in a matching is Aluffi torsion-free. Where, by a graph $G$ minus an edge $e$, we mean a graph resulting from $G$ which the edge $e$ is deleted and the vertices at the ends of $e$ are remaining.
    Note that, if $G$ is a complete graph minus a matching, then any vertex can be independent to at most only one other vertex. Therefore, item (i) of Theorem~\ref{graph} is not valid.
\item[\rm (iv)] The cycles $C_3$ and $C_4$ are Aluffi torsion-free. Because $C_3$ is a complete graph  and the  $C_4$ is the complete graph $K_4$ minus a maximal matching.
\item[\rm (v)] For each $n\geq 5$, the cycle $C_n$ is not Aluffi torsion-free. Let $n$ be even and $\{v_1,\ldots,v_n\}$ be the set of vertices of $G$ such that $v_i\sim v_{i+1}$ for $1\leq i\leq n-1$ and $v_n\sim v_1$.  Take $v_1$ and $v_2$ which are adjacent and $v_4, v_6, \ldots, v_{n-2}$ which are independent. Clearly condition (i) of Theorem~\ref{graph} is satisfied. Note that, $\hht(I(G))=\frac{n}{2}$ and degree of each vertex is 2. Moreover, $N(\{v_4, v_6, \ldots, v_{n-2}\})=\{v_3,v_5,\ldots, v_{n-1}\}$ which has cardinality $\frac{n}{2}-1$. This is condition (ii) of Theorem~\ref{graph}. \\
    If $n$ is odd, then, the vertices $v_1, v_2$ and $v_4, v_6, \ldots, v_{n-1}$ by the same argument as above, satisfy conditions of Theorem~\ref{graph}.
\item[\rm (vi)] Any path $P_n$ is not Aluffi torsion-free. It follows by the same argument as item (v) taking the same vertices.
\item[\rm (vii)] A star graph is not Aluffi torsion-free. Recall that a graph $G$ is called star if there is a vertex $v$, such that all other vertices are adjacent to $v$ and there is no any other edge. It is clear that $\hht(I(G))=1$. Then 1-minors of the Jacobian matrix of $I(G)$ is exactly the maximal ideal $\fm$. Hence, $I(G) = I(G) \cap \fm^2 \not\subseteq I(G)\fm$.
\end{itemize}
\end{Examples}

\begin{Remark}\rm
Let $G$ be a finite simple graph. Then, for $J=I(G)$, the edge ideal of $G$, Conjecture~\ref{conj1} holds.
\end{Remark}

\section{Other Examples}
In this section we consider some special class of varieties in algebraic geometry which have Aluffi torsion-free property. We assume that $k$ is algebraically closed of characteristic zero.
\begin{Examples}\rm
\begin{itemize}
  \item[\rm (i)]
Let $J\subset R=k[x,y,z]$ be the defining ideal of the monomial space curve with parametric equations $x=t^3,\ y=t^5,\ z=t^{7}$.
One knows that $J$ is a perfect ideal of codimension $2$ generated by  polynomials
$$F_1=x^4-yz, \ F_2=y^2-xz,\ F_3=x^3y-z^2$$
An easy calculation shows that the Jacobian ideal of $J$ is generated by monomials:
$$I=(x^4,x^3y,y^2,xz,yz,z^2)$$
We claim that the pair $J\subset I$ is Aluffi torsion-free. We show that the surjection in  Definition~\ref{torsionfree} is injective. By \cite[Lemma 2.1]{AA},
${\cal A}_{{R/J}}(I/J)\simeq \Rees_R(I)/(J,\tilde{J})\Rees_R(I)$. We compute the Rees algebra of $I$. It is
$$\Rees_R(I)=R[T_1,\ldots,T_6]/(I_1([\textbf{T}]\phi),\mathcal{J}),$$
where, $T=[T_1 \ T_2 \ \cdots \ T_6]$, $\phi$ is the first syzygy matrix of $I$,  $I_1([\textbf{T}]\phi)$ is the defining ideal of the symmetric algebra of $I$, and $\mathcal{J}$ is generated by the following polynomials.
$$T_5^2-T_3T_6,\ T_2T_4-T_1T_5,\ x^2T_4T_5-T_2T_6,\ x^2T_4^2-T_1T_6,\ x^2T_3T_4-T_2T_5,\ x^2T_1T_3-T_2^2$$
By the form of the generators of $I$ we have
$\tilde{J}=(T_1-T_5, T_2-T_6, T_3-T_4)$. Therefore
$${\cal A}_{{R/J}}(I/J)= R[T_1,\ldots,T_6]/(J,\tilde{J},I_1([\textbf{T}]\phi),\mathcal{J}).$$
In the other hand
$$\Rees_{R/J}(I/J)=R[T_1,\ldots,T_6]/(I_1([\textbf{T}]\widetilde{\phi}), T_5^2-T_4T_6,\ x^2T_4T_5-T_6^2,\ x^2T_4^2-T_5T_6),$$
where, $\widetilde{\phi}$ is the first syzygy matrix of the ideal $I/J$. Note that by definition, the defining ideal of the symmetric algebra of $I/J$ always
contains $\tilde{J}$. Thereof a direct calculation shows that the defining ideal of the Aluffi algebra of $I/J$ contains the defining ideal of the Rees algebra of $R/J$. This proves that the pair $J\subset I$ is Aluffi torsion-free.
\item[\rm (ii)] Let $J\subset R=K[x,y,z,w]$ be the defining ideal of the monomial space curve with parametric equations $x=t^3,\ y=t^4,\ z=t^{5}, w=t^{7}$. An easy calculation shows
$$J=(x^3-yz,y^2-xz,z^2-xw,x^2z-yw,xy-w).$$
While the Jacobian ideal of $J$ is
$$I=(xw, z^2, yz, xz, y^2,xy-w,x^3).$$
We use a computational argument as in Example (i) to show that the pair $J\subset I$ is Aluffi torsion-free. One has
$${\cal A}_{{R/J}}(I/J)= R[T_1,\ldots,T_7]/(J,\tilde{J},I_1([\textbf{T}]\phi),\mathcal{J}),$$
where, $\mathcal{J}=(T_2T_5-T_7^2, xT_5T_7-T_2^2, xT_5^2-T_2T_7)$, $\tilde{J}=(T_1-T_2, T_3-T_7, T_4-T_5, T_6)$ and $I_1([\textbf{T}]\phi)$ is the
defining ideal of the symmetric algebra of $I$. On the other hand, the defining ideal of the Rees Algebra of $I/J$ is generated by $J$, $I_1([\textbf{T}]\widetilde{\phi})$ and the following polynomials:
$$T_2T_5-T_7^2,\ xT_5T_7-T_2^2, \ wT_5^2-zT_7^2,\ xT_7^3-T_2^3,\ yzT_7^2-wT_2^2,$$$$xT_5^2-T_2T_7,\ yT_7^4-T_2^4,\ zT_7^5-T_2^5,\ wT_5T_7^5-T_6^2,\ wT_7^7-T_2^7$$
By a direct calculation, we can verify that the defining ideal of the Aluffi algebra contains the defining ideal of the Rees Algebra.
\end{itemize}
\end{Examples}
These examples motivate the following question.
\begin{Question}\rm
Let $J\subset k[x_1,\ldots,x_m] $ be defining ideal of the monomial space curve with parametric equations $x_1=t^{n_1},\ldots, x_m=t^{n_m}$,
 $\gcd(n_1,\ldots,n_m)=1$. Let $I$ be the Jacobian ideal of $J$. For which types of parametrization, the
pair $J\subset I$ is Aluffi torsion-free?
\end{Question}

\medskip

\begin{Remark}\rm
 Let $J$ be the defining ideal of a projective monomial curve in $\mathbb{P}_k^{d}$ defined
parametrically by $s^d,\ ts^{d-1},\ \ldots,\ t^{d-1}s,\ t^d$, which is the defining ideal of the rational normal curve. By Example~\ref{rational}(i), the pair $J\subset I$ is Aluffi torsion-free. We may derive the above question pretty generally.
Let $J$ be the defining ideal of a projective monomial curve in $\mathbb{P}_k^{n_d}$  parameterized by $s^{n_d},\ s^{n_d-n_1}t^{n_1},\ \ldots,\ t^{n_d-n_{d-1}}s^{n_{d-1}},\ t^{n_d}$ with positive integer $n_i$ such that $1\leq n_1<n_2<\ldots <n_d$. Let $I$ be the Jacobian ideal of $J$.
For which types of parametrizations, the pair $J\subset I$ is Aluffi torsion-free?
\end{Remark}

In \cite[Example 1.9.6]{AA} it is shown that if $J=(\partial f/\partial x_1, \ldots, \partial f/\partial x_n)$ is generated by partial derivations of $f=x_1\cdots x_n$, then the pair $J\subseteq I$ is Aluffi torsion-free where, $I$ is the Jacobian ideal of $J$.  This example stimulates us to consider the defining polynomial of a hyperplane arrangement and its Jacobian ideal.

A hyperplane arrangement $\mathcal{A}$ in $\mathbb{P}_k^{n-1}$ is a finite collection of hyperplanes. For each $H\in \mathcal{A}$, choose a linear polynomial $l_{H}$  as defining polynomial of it. Then, $Q=\displaystyle\prod_{H\in \mathcal{A}}l_{H}$ is the defining polynomial of the arrangement $\mathcal{A}$.

\begin{Examples}\rm
\begin{itemize}
  \item [\rm (i)] Let $\mathcal{A}$ be the hyperplane arrangement in $\mathbb{P}_k^{2}$ with equations
  $$x=y+z,\ y=x+z,\ z=x+y$$
  Setting $Q$ for the defining polynomial of $\mathcal{A}$. Then the Jacobian ideal of $Q$ is generated by three quadrics:
  $$
  \begin{array}{c}
    3x^2-2xy-y^2-2xz+2yz-z^2 \\
    x^2+2xy-3y^2-2xz+2yz+z^2 \\
    x^2-2xy+y^2+2xz+2yz-3z^2
  \end{array}
  $$
Note that $J$ is a codimension $2$ perfect ideal with syzygy matrix
$$M=
\begin{pmatrix}
  x+3y-3z & x-y+z \\
  3x+y-3z & x-y-z \\
  -2z & 2x-2y \\
  \end{pmatrix}
$$
A hard hand calculation shows that the Jacobian ideal $I$ of $J$ is $(x,y,z)^2$. Then by Corollary~\ref{HB}, the pair $J\subset I$ is Aluffi torsion-free.
\item [\rm (ii)] Let $\mathcal{A}$ be the hyperplane arrangement in $\mathbb{P}_k^{3}$ with equations
  $$x_1=x_2,\ x_2=x_3,\ x_3=x_4,\ x_4=x_1$$
  Let $Q$ be the defining polynomial of $\mathcal{A}$ and  $J$ be the Jacobian ideal of $Q$. Since $\hht(J)=2$, then  the Jacobian ideal  $I$ is generated by homogenous polynomials of degree $4$. To show that the pair $J\subset I$ is Aluffi torsion-free, we use \cite[Corollary 2.17]{AA}.  We need to calculate the relation type number of $I/J$.  It is $2$ by a direct machine computation. Hence, it is enough to  check
  that $J\cap I^2=JI$. An easy computation shows that this equality holds and the pair
  $J\subset I$ is Aluffi torsion-free. Note that the projective algebraic set defined by $J$ has one isolated singular point at $(1:1:1:1)$, then the Aluffi torsion-free property can be valid for some singular varieties too.
\end{itemize}
\end{Examples}
We wrap up with the following question.

\begin{Question}\rm
Let $Q=\displaystyle\prod_{H\in \mathcal{A}}l_{H}$ be the defining polynomial of the hyperplane arrangement $\mathcal{A}$.
Let $J$ be the Jacobian ideal of $Q$. If $I$ is the Jacobian ideal of $J$, when is $J\subset I$ an Aluffi torsion-free pair.
\end{Question}

\end{document}